\documentclass{amsart}

\usepackage{amssymb, latexsym,pdfsync}
\theoremstyle{plain}
\newtheorem{theorem}{Theorem}
\newtheorem{corollary}{Corollary}
\newtheorem{lemma}{Lemma}

\theoremstyle{definition}

\newtheorem{remark}{Remark}
\newtheorem{example}{Example}

\newcommand{\FF}{\mathbb{F}}

\newcommand{\Fq}{\mathbb{F}_q}

\newcommand{\D}{\mathcal D}

\def\Fq{{\mathbb{F}}_q}

\def\End{\mathrm{End}}

\def\PG{\mathrm{PG}}

\def\dim{\mathrm{dim}}

\def\cP{\mathcal{P}}

\def\Tr{\mathrm{Tr}}

\begin{document}

\title{Dimensional dual hyperovals in classical polar spaces}
\author{John Sheekey}
\maketitle

\begin{abstract}
In this paper we show that $n$-dimensional dual hyperovals cannot exist in all but one classical polar space of rank $n$ if $n$ is \emph{even}. This resolves a question posed by Yoshiara.
\end{abstract}

\section{Definitions and preliminaries}

An \emph{$n$-dimensional dual arc} $\D$ in a vector space $V(N,q)$ over a finite field $\Fq$ is a set of $n$-dimensional subspaces such that
\begin{enumerate}
\item
each two intersect in exactly a one-dimensional space;
\item
no three intersect non-trivially.
\end{enumerate}
It is clear that $|\D| \leq \frac{q^n-1}{q-1}+1$. For let $S$ be any element of $\D$. Then the other elements of $\D$ intersect $S$ in distinct one-dimensional subspaces, of which there are $\frac{q^n-1}{q-1}$. If $\D$ meets this bound, it is called an \emph{$n$-dimensional dual hyperoval}. We will sometimes use the shorthand $n$-DA and $n$-DHO.

For background and a recent survey of known results and applications, we refer to \cite{Yoshiara2006}. Note that the definition therein are in terms of projective spaces, but here we use vector space terminology and notation, following \cite{DeKa2015}. In this paper we will mostly consider the case $N=2n$. In \cite{DeKa2015}, this is required in the definition, but we will not impose this restriction here.

It is known that $n$-dimensional dual hyperovals exist in $V(2n,q)$ for all $n$ and all $q$ even, see for example \cite{Yoshiara2006}. It is an open problem whether any can exist when $q$ is odd.

In this paper we will consider $n$-dimensional dual arcs in \emph{polar spaces}, that is, where $\D$ consists of maximum totally isotropic subspaces with respect to some nondegenerate form on $V(N,q)$. Necessarily then we have that $N \in \{2n,2n+1,2n+2\}$. 

It is known \cite{Yoshiara2006} that there exist $n$-dimensional dual hyperovals in the hyperbolic quadric $Q^{+}(2n-1,q)$ for all $n$ odd and $q=2$ (see Section \ref{sec:polar} for notation). Furthermore, there exists a $3$-dimensional dual hyperoval in the hermitian variety $H(5,4)$, the \emph{Mathieu dual hyperoval}. 

In \cite{Yoshiara2006} Problem 4.7, the following (paraphrased) question is asked.

{\it Does the existence of an $n$-dimensional dual hyperoval in a polar space imply that $n$ is odd?}

Taniguchi \cite{Taniguchi2009} proved that $n$-dimensional ``alternating doubly dual hyperovals'' exist in $V(2n,2)$ if and only if $n$ is odd. Dempwolff \cite{Dempwolff2015}, showed that $n$-dimensional ``symmetric doubly dual hyperovals'' exist only if $n$ is odd. We will see in Section \ref{sec:demp} that the existence of such implies the existence of an $n$-dimensional dual hyperovals in the symplectic space $W(2n-1,q)$.

We respond now to these questions with the following theorem.

\begin{theorem}\label{thm:main}
Suppose $\D$ is an $n$-dimensional dual hyperoval in a polar space $\cP$ of rank $n$. Then either $n$ is odd or $\cP$ is an elliptic quadric. 
\end{theorem}

The result is a simple application of a theorem of Vanhove.

\section{Polar spaces}
\label{sec:polar}

A classical \emph{polar space} $\cP$ is the geometry of totally singular subspaces with respect to some non-degenerate quadratic form on $V(N,q)$, or totally isotropic with respect to some non-degenerate symplectic or sesquilinear form on $V(N,q)$. The \emph{rank} of $\cP$ is the maximum (vector space) dimension of a subspace of $\cP$. If every $(n-1)$-dimensional space of a polar space of rank $n$ is contained in precisely $q^e+1$  totally isotropic $n$-dimensional spaces, then $\cP$ is said to have \emph{parameters} $(q,q^e)$. See for example \cite{DeBKlMe2011} for background. We tabulate the relevant polar spaces of rank $n$ here. 

\begin{center}
\begin{tabular}{|c|c|c|c|c|c|}
\hline
Name & form & Notation & Ambient vector & Parameters & $e$\\ 
&&& space &&\\
\hline
Hyperbolic quadric & quadratic & $Q^{+}(2n-1,q)$ & $V(2n,q)$ & $(q,1)$& $0$\\
\hline
Parabolic quadric & quadratic & $Q(2n,q)$ & $V(2n+1,q)$ & $(q,q)$& $1$\\
\hline
Elliptic quadric & quadratic & $Q^{-}(2n+1,q)$ & $V(2n+2,q)$ & $(q,q^2)$& $2$\\
\hline
Symplectic space & symplectic & $W(2n-1,q)$ & $V(2n,q)$ & $(q,q)$& $1$\\
\hline
Hermitian variety & sesquilinear & $H(2n-1,q^2)$ & $V(2n,q^2)$ & $(q^2,q)$& $1/2$\\
\hline
Hermitian variety & sesquilinear & $H(2n,q^2)$ & $V(2n+1,q^2)$ & $(q^2,q^3)$& $3/2$\\
\hline
\end{tabular}
\end{center}
Note that $n$-dimensional dual hyperovals in polar spaces defined by a quadratic form are often referred to as being of {\it orthogonal type}.

If $q$ is even, $W(2n-1,q)$ is isomorphic to $Q(2n,q)$, and contains $Q^{+}(2n-1,q)$. 

\begin{example}
Yoshiara defined in \cite{Yoshiara1999} the following $n$-dimensional dual hyperovals in $V(2n,2)$, and showed in \cite{Yoshiara2005} that they lie in $Q^{+}(2n-1,2)$ (and hence $W(2n-1,2)$) if and only if $n$ is \emph{odd}. Let $h$ be an integer coprime to $n$, and for each $t\in \FF_{2^n}$ define
\[
S_t = \{(x,x^{2^{-2h}}t+xt^{2^h}):x \in \FF_{2^n}\}.
\]
Then $\D := \{S_t:t \in \FF_{2^n}\}$ is an $n$-dimensional dual hyperoval in $Q^+(2n-1,2)$ (and $W(2n-1,2)$), where the quadratic form on $V(2n,2)$ is
\[
(a,b) \mapsto \Tr(ab^{2^h}),
\]
and the associated symmetric (alternating) bilinear form on $V(2n,2)$ is 
\[
((a,b),(c,d)) \mapsto \Tr(ad^{2^h} - bc^{2^h}).
\]
Dempwolff and Kantor \cite{DeKa2015} gave a geometric construction leading to many inequivalent examples in $Q^{+}(2n-1,2)$. Dempwolff \cite{Dempwolff2013} gave further examples in $W(2n-1,2)$ which cannot lie in  $Q^+(2n-1,2)$.
\end{example}

\begin{example}
There exists a $3$-dimensional dual hyperoval in $V(6,4)$ which lies in the polar space $H(5,4)$ known as the \emph{Mathieu dual hyperoval}, see e.g. \cite{DelFra2000}.
\end{example}

To the author's knowledge, no examples in other polar spaces are known. Del Fra \cite{DelFra2000} showed that the only $3$-dimensional dual hyperovals in a polar space are the above examples. 

Yoshiara \cite{Yoshiara2005} showed that $n$-dimensional dual hyperovals can exist in $Q^{+}(2n-1,q)$ only if $n$ is odd.

\section{Dual polar graphs and Main result}

Given a polar space $\cP$ of rank $n$, we define the \emph{dual polar graph} $\Gamma_{\cP}$, whose vertices are the $n$-spaces of $\cP$, and where two vertices are adjacent if their intersection has dimension $n-1$. Many properties of this graph are know, see for example \cite{BrCoNe}, \cite{VanhoveThesis}. 

For a set $\D$ of $n$-spaces of $\cP$, the \emph{inner distribution} is an $(n+1)$-tuple of integers $a = (a_0,a_1,\ldots,a_n)$, where
\[
a_i = \frac{\{(S,T):S,T \in \D~\mid~\dim(S \cap T) = n-i\}}{|\D|}.
\]
Equivalently, if we view $\D$ as a subset of $\Gamma_{\cP}$, and let $d(S,T)$ denote the distance function on $\Gamma$, then
\[
a_i = \frac{\{(S,T):S,T \in \D~\mid~d(S,T) = i\}}{|\D|}.
\]

In \cite[Lemma 3.2]{Van10}, the following was proved.
\begin{theorem}[Vanhove]
Let $\cP$ be a classical polar space of rank $n$ with parameters $(q,q^e)$, and let $\D$ be a set of $n$-spaces in $\cP$ with inner distribution $(a_0,a_1,\ldots,a_n)$. Then
\[
\sum_{i=0}^n \left(-\frac{1}{q^e}\right)^i a_i \geq 0.
\]
\end{theorem}

Now suppose $\D$ is a dimensional dual arc in $\cP$. Then it is clear that
\begin{align*}
a_0 &= 1\\
a_{n-1} &= |\D|-1\\
a_i &= 0 \textrm{ otherwise}.
\end{align*}
Hence we get that
\[
1 + \left(-\frac{1}{q^e}\right)^{n-1}(|\D|-1) \geq 0,
\]
and so if $n$ is even,
\[
|\D| \leq q^{(n-1)e} + 1.
\]
Hence we get an upper bound for an $n$-dimensional dual arc in each classical polar space.
\begin{theorem}
\label{thm:bound}
Suppose $\D$ is an $n$-dimensional dual arc in $\cP$, and suppose $n$ is even. Then we have the following upper bounds on $|\D|$.

\begin{tabular}{|c|c|c|c|c|c|}
\hline
$\cP$ & Ambient vector & Parameters & $e$ & $|\D|\leq$ & Size of DHO\\ 
& space &&&&\\
\hline
 $Q^{+}(2n-1,q)$ & $V(2n,q)$ & $(q,1)$& $0$ & $2$ & $\frac{q^n-1}{q-1}+1$\\
\hline
 $Q(2n,q)$ & $V(2n+1,q)$ & $(q,q)$& $1$&$q^{n-1}+1$ & $\frac{q^n-1}{q-1}+1$\\
\hline
 $Q^{-}(2n+1,q)$ & $V(2n+2,q)$ & $(q,q^2)$& $2$&$q^{2n-2}+1$ & $\frac{q^n-1}{q-1}+1$\\
\hline
 $W(2n-1,q)$ & $V(2n,q)$ & $(q,q)$& $1$&$q^{n-1}+1$ & $\frac{q^n-1}{q-1}+1$\\
\hline
 $H(2n-1,q^2)$ & $V(2n,q^2)$ & $(q^2,q)$& $1/2$&$q^{n-1}+1$ & $\frac{q^{2n}-1}{q^2-1}+1$\\
\hline
 $H(2n,q^2)$ & $V(2n+1,q^2)$ & $(q^2,q^3)$& $3/2$&$q^{3(n-1)/2}+1$ & $\frac{q^{2n}-1}{q^2-1}+1$\\
\hline
\end{tabular}
\end{theorem}

{\it Proof of Theorem \ref{thm:main}:} This now now follows immediately by comparing  the above upper bounds on $n$-dimensional dual arcs (fourth column) with the required size of an $n$-dimensional dual hyperoval (fifth column).

\begin{remark}
Note that an $n$-dimensional dual hyperoval is a special case of a \emph{constant-distance, constant-dimension subspace code} \cite{GoLaShVa2014}, \cite{Ihringer2014}, or equivalently, a clique in the graph $\Gamma_{n-1}$, where $\Gamma_i$ is the graph whose vertices are the vertices of $\Gamma$, and whose edges are between vertices at distance $i$ in $\Gamma$. Note however that not every clique of the correct size in $\Gamma_{n-1}$ gives rise to an $n$-dimensional dual hyperoval. As the proof of Theorem \ref{thm:bound} does not use the fact that no three spaces intersect nontrivially, the same bounds hold for the relevant constant-distance subspace codes in each polar spaces.

This is the same method used by Vanhove in \cite{Van09} to prove that the maximum size of a partial spread in $H(2n-1,q)$ is $q^n+1$ if $n$ is odd.
\end{remark}

\begin{remark}
This table does not imply any results for dimensional dual hyperovals in elliptic quadrics $Q^{-}(2n+1,q)$. This problem seems to require a different approach. Note that such objects do not satisfy the definition of a dimensional dual hyperoval in \cite{DeKa2015}.
\end{remark}

\section{Alternating and symmetric doubly dual hyperovals}
\label{sec:demp}
 
An $n$-dimensional dual hyperoval $\D$ in $V(2n,q)$ is said to be ``doubly dual'' if $\D^{\perp} := \{S^{\perp}:S \in \D\}$ is also an $n$-dimensional dual hyperoval, where $\perp$ is some nondegenerate polarity. Note that if $\D$ lies in some polar space $\cP$, it is doubly dual: we take $\perp$ to be the polarity defined by the quadratic or sesquilinear form associated to $\cP$, whence $S^{\perp}=S$ for all maximum subspaces $S$ in $\cP$. However, the converse is not necessarily true.

In \cite{Dempwolff2015} the concept of a (bilinear) symmetric doubly dual hyperoval was introduced, and it was proved that such objects can not exist in $V(2n,q)$ for $n$ even. We will now show that the existence of this implies the existence of an $n$-dimensional dual hyperoval in symplectic polar space.

Suppose there is some injective linear map $\beta:V(n,q)\rightarrow \End(V(n,q))$. Let us represent elements of $V(2n,q)$ with elements of $V(n,q)\times V(n,q)$. For each $y\in V(n,q)$, define an $n$-dimensional subspace $S_{y} = \{(x,\beta(y)(x)):x \in V(n,q)\}$ of $V(2n,q)$, and define $\D_{\beta} = \{S_y:y \in V(n,q)\}$. If $\D_{\beta}$ is an $n$-dimensional dual hyperoval, then it is called a \emph{bilinear dual hyperoval}. Note that this can occur only if $q=2$. 

Define $\beta^o:V(n,q)\rightarrow \End(V(n,q))$ by $\beta^o(x)(y) = \beta(y)(x)$. If $\beta=\beta^o$, that is if $\beta(y)(x) = \beta(x)(y)$ for all $x,y \in V(n,q)$, then $\D_{\beta}$ is called a \emph{symmetric dual hyperoval}. If furthermore $\beta(x)(x)=0$ for all $x$, then $\D_{\beta}$ is called an \emph{alternating dual hyperoval}.

Let $\langle,\rangle : V(n,q)\times V(n,q) \rightarrow \Fq$ be a nondegenerate symmetric bilinear form on $V(n,q)$. Let $t$ denote the adjoint operator with respect to this form, i.e. $\langle x,f(y)\rangle = \langle f^t(x),y\rangle$ for all $x,y \in V(n,q)$, and define $\beta^t:V(n,q)\rightarrow \End(V(n,q))$ by $\beta^t(x) = \beta(x)^t$.

Taniguchi \cite{Taniguchi2009} showed that alternating doubly dual hyperovals exist in $V(2n,2)$ if and only if $n$ is odd. Dempwolff \cite{Dempwolff2015} improved this by showing that symmetric doubly dual hyperovals exist in $V(2n,2)$ only if $n$ is odd. We now show that this result follows also from Theorem \ref{thm:main} of this paper.

The following lemma follows from \cite{Edel2010}, and from \cite[Lemma 3.8]{Dempwolff2013}.

\begin{lemma}
If there exists a symmetric bilinear $n$-dimensional doubly dual hyperoval $\D$ in $V(2n,2)$, then there exists an $n$-dimensional dual hyperoval in $W(2n-1,2)$.
\end{lemma}

Combining this with Theorem \ref{thm:main} immediately gives us the following corollary.

\begin{corollary}
There exists a symmetric bilinear $n$-dimensional doubly dual hyperoval $\D$ in $V(2n,2)$ only if $n$ is odd.
\end{corollary}
%
%
%
%

Note that Theorem \ref{thm:main}, applied to $W(2n-1,q)$, does not require either bilinearity or that $q=2$, and so this result is more general than the results of Taniguchi and Dempwolff.

Dempwolff further conjectured in \cite{Dempwolff2013} that $n$-dimensional doubly dual hyperovals over $\FF_2$ exist only if $n$ is odd. This remains an open problem.

%
%

\section{Acknowledgements}

The results of this paper were developed in discussions with Fr\`ed\`eric Vanhove prior to his tragic early passing. The author is heavily indebted to Fr\`ed\`eric for this work, and this paper is dedicated to his memory.

The author is supported by the Fund for Scientific Research Flanders (FWO -- Vlaanderen).


\begin{thebibliography}{9}

\bibitem{BrCoNe}
Brouwer, A.E., Cohen, A.M., Neumaier, A.: {\it Distance regular graphs}, 
 Springer-Verlag, New York, 1989.

\bibitem{DeBKlMe2011}
De Beule, J.; Klein, A.; Metsch. K.: Substructures of finite classical polar spaces. Chapter in {\it Current research topics in Galois Geometry} (Editors J. De Beule and L. Storme), NOVA Academic Publishers, New York, 2011.

\bibitem{Dempwolff2013}
Dempwolff, U.: Dimensional doubly dual hyperovals and bent functions, {\it Innovations in Incidence Geometry} 13 (2013) 149-178.
 
\bibitem{Dempwolff2015}
Dempwolff, U: Symmetric doubly dual hyperovals have an odd rank, {\it Des. Codes Cryptogr.} 74 (2015) 153-157.

\bibitem{DeKa2015}
Dempwolff, U.; Kantor, W.M.: Orthogonal dual hyperovals, symplectic spreads and orthogonal spreads, {\it J. Alg. Comb.} 41(2015) 83-108.

\bibitem{DelFra2000}
Del Fra, A: On d-dimensional dual hyperovals, {\it Geom. Dedicata} 79 (2000), 157-178.

\bibitem{Edel2010}
Edel, Y.: On some representations of quadratic APN functions and dimensional dual hyperovals, {\it RIMS Kokyuroku} 1687 (2010) 118-130.

\bibitem{GoLaShVa2014}
Gow, R.; Lavrauw, M.; Sheekey, J.; Vanhove, F.: Constant rank-distance sets of Hermitian matrices and partial spreads in Hermitian polar spaces, {\it Electron. J. Combin.} 21 (2014) Paper 1.26, 19 pp. 

\bibitem{Ihringer2014}
Ihringer, F.: A new upper bound for constant distance codes of generators on Hermitian polar spaces of type $H(2d-1,q^2)$, {\it J. Geom.} 105 (2014) 457-464.

\bibitem{Taniguchi2009}
Taniguchi, H.: On the duals of certain d-dimensional dual hyperovals in $\PG(2d+1,2)$, {\it Finite Fields Appl.} 15 (2009) 673-681.

\bibitem{Van09}
Vanhove, F.: The maximum size of a partial spread in $H(4n + 1, q^2)$ is $q^{2n+1} + 1$, {\it Electron.
J. Combin.}, 16 (2009), 1--6.

\bibitem{VanhoveThesis}
Vanhove, F.: {\it Incidence geometry from an algebraic graph theory point of view}, Ph.D. Thesis (2011).

\bibitem{Van10}
Vanhove, F.: Antidesigns and regularity of partial spreads in dual polar graphs, {\it J. Combin. Des.} 19 (2011), 202-216.


\bibitem{Yoshiara1999}
Yoshiara, S.: A family of d-dimensional dual hyperovals in $\PG(2d + 1, 2)$, {\it Europ. J. Combin.} 20 (1999), 589-603.

\bibitem{Yoshiara2006}
Yoshiara, S.: Dimensional dual arcs: a survey. {\it Finite geometries, groups, and computation}, 247-266, Walter de Gruyter GmbH \& Co. KG, Berlin, 2006.

\bibitem{Yoshiara2005}
Yoshiara, S.: Some remarks on dimensional dual hyperovals of polar type, {\it Bull. Belg. Math. Soc. Simon Stevin} 12 (2005) 925-939.

\end{thebibliography}
\end{document}